\documentstyle{amsppt}
\magnification=\magstep1
\pagewidth{6truein}
\pageheight{8.5truein}
\hcorrection{5mm}
\vcorrection{12mm}
\baselineskip=14pt
\parindent=12mm

\def\BC{\Bbb C}

\def\Def{Definition}

\def\nbd{neighborhood}

\def\a{\alpha}

\def\op{\operatorname}

\topmatter
\title A proof of the proportionality theorem
\endtitle

\author  J.-P. Brasselet, J. Seade and T. Suwa\endauthor
\address Institut de Math\'ematiques de Luminy, UPR 9016 CNRS, Campus
de Luminy - Case 930, 13288 Marseille Cedex 9, France
\endaddress
\email jpb$\@$iml.univ-mrs.fr\endemail
\address
Instituto de Matem\'aticas, UNAM, Unidad Cuernavaca,
Av. Universidad s/n, Colonia Lomas de Chamilpa,
C.P. 62210, Cuernavaca, Morelos,
M\'exico
\endaddress
\email jseade$\@$matem.unam.mx
\endemail
\address Department of Information Engineering, Niigata University, 
2-8050,\newline Ikarashi, Niigata 950-2181,
Japan\endaddress
\email suwa$\@$ie.niigata-u.ac.jp\endemail
\date
November 22, 2005\enddate

\thanks Research partially supported by  CONACYT and DGAPA-UNAM    (JS) 
and by JSPS  (JPB and TS). 
\endthanks

\endtopmatter
\document
\NoRunningHeads
The Proportionality Theorem of [BS] establishes a deep relation between the 
Schwartz index of stratified vector fields on Whitney stratified complex analytic varieties and 
the liftings of these vector fields to sections of the Nash bundle. This is one of the key ingredients
for the proof in [BS] that the Alexander isomorphism carries the Schwartz classes of 
singular varieties to the corresponding MacPherson classes. 

Recently we considered the corresponding problem for $1$-forms and gave a proof of the 
Proportionality Theorem in [BSS]. It is inspired by the original proof for vector fields in [BS],
however is shorter and more direct. Then, coming back to the case of vector fields, we realized that 
the original proof, for vector fields and for frames as well, can be substantially simplified.

The purpose of this note is to give a direct and self-contained proof of the Proportionality Theorem
in order to facilitate the understanding of this important theorem.

Part of this work was done during  the authors' stay at ICTP, Trieste Italy, and first author's
stay at Niigata University. We thank these institutions for their financial support and hospitality.

\bigpagebreak

\noindent
{\bf 1. The case of vector fields}
\medpagebreak

Let $V$ be a {complex} 
analytic variety of pure dimension $n$ in a complex manifold $M$ of 
dimension
$m$. We endow $M$ with a Whitney stratification $\{V_\a\}$ adapted to $V$.  
Let $SV$ denote
the union of the tangent bundles of all the strata; $SV=\bigcup_\a TV_\a$, 
which is thought
of as a subset of the tangent bundle $TM$ of $M$. A {\it stratified vector 
field} is a section of
$TM$ whose image is in $SV$. Here $TM$ denotes the complex tangent bundle of 
$M$, however,
recall that it is canonically isomorphic to the real tangent bundle as a 
real bundle.

The Nash modification  and the Nash bundle are constructed as follows. 
Let $\op{Sing}(V)$ denote the singular set and $V_{\op{reg}}=V\setminus 
\op{Sing}(V)$
the regular part of $V$. We have a map $\sigma:V_{\op{reg}}\to G_n(TM)$ into the 
Grassmann bundle
of {complex} 
$n$ planes in $TM$, which assigns to each point $p$ in $V_{\op{reg}}$
the subspace $T_pV_{\op{reg}}$ of $T_pM$.
The Nash modification $\widetilde V$ is the closure of the image of $\sigma$ together with
 the restriction $\nu:\widetilde V\to V$ of the 
projection $G_n(TM)\to M$ to $\widetilde V$.
The Nash bundle $\widetilde T\to \widetilde V$  is the 
restriction
to $\widetilde V$ of the tautological bundle over $G_n(TM)$. Note that 
$\widetilde T$ is a subbundle
of $\nu^*TM|_V$ and is isomorphic to 
$TV_{\op{reg}}$ away from $\nu^{-1}(\op{Sing}(V))$.

If $v$ is a non-vanishing
stratified vector field on a subset $A$ in $V$, then by the Whitney 
condition (a), 
it can be lifted to a non-vanishing section $\widetilde v$ of $\widetilde T$ 
over $\nu^{-1}A$ {(see [BS])}. Namely, the 
pull-back $\nu^*v$, which is a priori a section of $\nu^*TM|_A$, is in fact 
a section of
$\widetilde T$. {We denote it}  by $\widetilde v$ to emphasize that it 
is a section of $\widetilde T$.

We take a  closed ball 
$\Bbb B$ in $M$ around a point $p$ in $V$ and let $\Bbb S
=\partial \Bbb B$. If $v$ is a stratified vector field, which is non-vanishing on
$\Bbb S\cap V$, 
it can be lifted to a non-vanishing section $\widetilde v$
of $\widetilde T$ over 
$\nu^{-1}(\Bbb S\cap V)$, as noted above. Let $o(\widetilde v)$ denote the 
class in 
$H^{2n}(\nu^{-1}(\Bbb B\cap V),\nu^{-1}(\Bbb S\cap V))$
of the obstruction cocycle to extending this to a non-vanishing section over 
$\nu^{-1}(\Bbb B\cap V)$.

\definition{Definition 1.1} The {\it local Euler obstruction} 
$\op{Eu}(v,V;p)$
of a stratified vector field $v$ at an isolated
singularity  $p$ is the integer obtained by evaluating $o(\widetilde v)$ on 
the orientation cycle
$[\nu^{-1}(\Bbb B\cap V),\nu^{-1}(\Bbb S\cap V)]$, for sufficiently small $\Bbb B$.
\enddefinition

As a particular vector field  satisfying the previous  situation,
one can consider a radial vector field, $v_{\op{rad}}$, i.e. a vector field 
pointing outwards  $\Bbb B$ along $\Bbb S$.
By [BS],  the local Euler obstruction $\op{Eu}(v_{\op{rad}},V;p)$ 
for a {stratified} radial vector
field  coincides with the {\it local Euler obstruction of $V$ at $p$}, 
$\op{Eu}(V,p)$, introduced by
MacPherson [M]; we take this as our definition of  $\op{Eu}(V,p)$.

\medpagebreak

Now let $v_\a$ be a vector field on a stratum $V_\a$ with an isolated singularity 
at $p\in V_\a$.
Recall that by the radial extension process of M.-H. Schwartz [Sc], $v_\a$ 
can be extended 
to a  vector field $v'_\a$ {defined} in a \nbd\ of $p$ in $M$.
In the rest of this section, we denote the vector field $v'_\a$ simply by $v$.
The vector field $v$ is stratified with an isolated 
singularity at $p$. Moreover, it has the property that its Poincar\'e-Hopf
index $\op{Ind}_{\op{PH}}(v,M;p)$ in $M$ coincides with the Poincar\'e-Hopf
index $\op{Ind}_{\op{PH}}(v_\a,V_\a;p)$ of $v_\a$ in $V_\a$, provided that $\dim 
V_\a>0$.
If $\dim V_\a=0$, then $v$ is a radial vector field in the usual sense, and  
we have 
$\op{Ind}_{\op{PH}}(v,M;p)=1$.

\definition{Definition 1.2}  Let  $v$ be as above. 
The {\it Schwartz index} $\op{Ind}_{\op{Sch}}(v,V;p)$ of $v$ relative to $V$ at $p$ is the 
Poincar\'e-Hopf
index $\op{Ind}_{\op{PH}}(v,M;p)$ of $v$ as a  vector field on $M$, or equivalently, 
if $\dim V_\a>0$,  the Poincar\'e-Hopf
index $\op{Ind}_{\op{PH}}(v_\a,V_\a;p)$ of $v_\a$ on $V_\a$.
\enddefinition

First we give the proposed proof of the Proportionality Theorem in [BS] in the case of
vector fields. 

\proclaim{Theorem 1.3} Let $V_\a\subset V$ be a stratum and $v_\a$ a vector field on 
$V_\a$ with an isolated 
singularity at $p$. Let $v$ denote the radial extension of $v_\a$. Then we have
$$
\op{Eu}(v,V;p)=\op{Eu}(V,p)\cdot \op{Ind}_{\op{Sch}}(v,V;p).
$$
\endproclaim

\demo{Proof} In the sequel, we denote by $T^\times M$ and $\widetilde T^\times$ the 
spaces obtained from
$T M$ and $\widetilde T$ by removing the zero sections. 

First  we recall (\Def\ 1.2) that
$$
\op{Ind}_{\op{Sch}}(v,V;p)=\op{Ind}_{\op{PH}}(v, M;p).\tag1.4
$$

Let $v_{\op{rad}}$ 
denote a stratified radial vector field at $p$. Then, by definition of 
$\op{Ind}_{\op{PH}}(v,M;p)$, there is a homotopy
$$
\Psi:\Bbb S\times [0,1]@>>>T^\times M|_{\Bbb S}
$$
such that
$$
\partial \op{Im}\Psi
=v(\Bbb S)-\op{Ind}_{\op{PH}}(v,M;p)\cdot v_{\op{rad}}(\Bbb S)
$$
as chains in  $T^\times  M|_{\Bbb S}$. Moreover, since the vector fields
$v$ and  $v_{\op{rad}}$ are both stratified and the stratification is Whitney,
 we can assume that the above homotopy $\Psi$ is through stratified vector fields, i.e., 
we may choose $\Psi$ 
so that 
$$
\op{Im}\Psi\subset SV.  \tag1.5
$$

The restriction of $\Psi$ gives a homotopy
$$
\psi:(\Bbb S\cap V)\times [0,1]@>>>T^\times  M|_{\Bbb S\cap V}
$$
such that (cf. (1.4))
$$
\partial \op{Im}\psi
=v(\Bbb S\cap V)-\op{Ind}_{\op{Sch}}(v,V;p)\cdot 
v_{\op{rad}}(\Bbb S\cap V).
$$

We can lift  $v$ and $v_{\op{rad}}$ to sections 
$\nu^* v$ and $\nu^* v_{\op{rad}}$ of $\nu^*T^\times M|_{\nu^{-1}(\Bbb S\cap V)}$;
we  can also
lift $\psi$ to a homotopy
$$
\nu^*\psi:\nu^{-1}(\Bbb S\cap V)\times [0,1]
@>>>\nu^*T^\times  M|_{\nu^{-1}( \Bbb S\cap V)}
$$
and we have
$$
\partial \op{Im}\nu^*\psi
=\nu^*v(\nu^{-1}(\Bbb S\cap V))-\op{Ind}_{\op{Sch}}(v,V;p)\cdot 
\nu^*v_{\op{rad}}(\nu^{-1}(\Bbb S\cap V))
$$
as chains in $\nu^*T^\times M|_{\nu^{-1}(\Bbb S\cap V)}$. By (1.5), everything 
can be restricted to $\widetilde T\subset \nu^*TM$ to get a homotopy
$$
\widetilde\psi:\nu^{-1}(\Bbb S\cap V)\times [0,1]
@>>>\widetilde T^\times |_{\nu^{-1}(\Bbb S\cap V)}
$$
and we  have
$$
\partial\op{Im}\widetilde\psi=\widetilde v(\nu^{-1}(\Bbb S\cap V))-
\op{Ind}_{\op{Sch}}(v,V;p)
\cdot\widetilde v_{\op{rad}}(\nu^{-1}(\Bbb S\cap V))
$$
as chains in $\widetilde T^\times |_{\nu^{-1}(\Bbb S\cap V)}$. 

Taking a triangulation or a cellular decomposition of 
$\nu^{-1}(\Bbb B\cap V)$ and extending the homotopy $\tilde\psi$ to the $(2n-1)$-skeleton
of the decomposition, we see that the obstruction to extending $\tilde v$ is 
$\op{Ind}_{\op{Sch}}(v,V;p)$ times the obstruction to extending $\widetilde v_{\op{rad}}$.
By definition of the Euler obstructions, we have 
the theorem.\qed
\enddemo

\bigpagebreak

\noindent
{\bf 2. The case of frames}
\medpagebreak

Let $M$, $V$ and $\{V_\a\}$ be as in Section 1. We take a triangulation $(K)$ of $M$
compatible with
the stratification and let $(D)$ denote the cellular decomposition dual to $(K)$. Note that 
the cells
in $(D)$ are transverse to $V$ and $V_\a$ so that if $\sigma$ denotes a cell of real dimension
$2s$, then $\sigma\cap V$ and $\sigma\cap V_a$ are of dimensions $2(s-m+n)$ and 
$2(s-m+n_\a)$,
respectively, where $n_\a=\dim_{\BC}V_\a$.

In the sequel, an {\it $r$-field}  means a collection 
$v^{(r)}=(v_1,\dots,v_r)$ of $r$ vector fields. A singular
point of $v^{(r)}$ is a point where the vectors fail to be linearly independent
(over the complex numbers).
An {\it $r$-frame} is an $r$-field without singularity. 

Let  $\nu:\widetilde V\to V$ be
the Nash modification of $V$. 
Let $\sigma$ be a cell of dimension $2(m-r+1)$ and  $v^{(r)}$ a stratified $r$-field on 
$\sigma\cap V$ with an isolated singularity at the barycenter $p$ of $\sigma$.
Since $v^{(r)}$ is non-singular on $\partial\sigma\cap V$, 
it can be lifted to an $r$-frame  
$\widetilde v^{(r)}$ of $\widetilde T$ over 
$\nu^{-1}(\partial\sigma\cap V)$, as in the case of vector fields.  
Let $o(\widetilde v^{(r)})$ denote the class in 
$H^{2(n-r+1)}(\nu^{-1}(\sigma\cap V),\nu^{-1}(\partial\sigma\cap V))$
of the obstruction cocycle to extending this to an $r$-frame over 
$\nu^{-1}(\sigma\cap V)$.

\definition{Definition 2.1} The {\it local Euler obstruction} $\op{Eu}(v^{(r)},V;p)$
of a stratified $r$-field $v^{(r)}$ at an isolated
singularity  $p$ is the integer obtained by evaluationg $o(\widetilde v^{(r)})$ on the orientation cycle
$[\nu^{-1}(\sigma\cap V),\nu^{-1}(\partial\sigma\cap V)]$.
\enddefinition

\medpagebreak

Now let $v^{(r)}_\a$ be an $r$-field on $\sigma\cap V_\a$ with an isolated 
singularity
at the barycenter $p$ of $\sigma$, where $\sigma$ is a cell of dimension
$2(m-r+1)$. We may construct an $r$-field ${v_\a^{(r)}}'$ on $\sigma$ by the radial extension
process of Schwartz [Sc]. In the sequel, we denote ${v_\a^{(r)}}'$ simply by
$v^{(r)}$. The $r$-field $v^{(r)}$ is stratified with an isolated singularity
at $p$. Moreover, it has the property that its index $\op{Ind}(v^{(r)},M;p)$ (the
obstruction to extending $v^{(r)}$ as an $r$-frame of $TM$ on $\sigma$) coincides with
the index $\op{Ind}(v_\a^{(r)},V_\a;p)$ (the
obstruction to extending $v_\a^{(r)}$ as an $r$-frame of $TV_\a$ on $\sigma\cap V_\a$), provided
that $\dim_{\BC} V_\a \ge r$.
If $\dim V_\a=r-1$, then $v^{(r)}$ is a radial $r$-field in the usual sense, and  
we have 
$\op{Ind}_{\op{PH}}(v^{(r)},M;p)=1$.

\definition{Definition 2.2} Let $v^{(r)}$ and $\sigma$ be as above. The {\it Schwartz index} 
$\op{Ind}_{\op{Sch}}(v^{(r)},V;p)$ of $v^{(r)}$ relative to $V$
at an isolated singularity $p$ is the
index $\op{Ind}(v^{(r)},M;p)$ of $v^{(r)}$ as an $r$-field on $\sigma$, or equivalently,  the 
index $\op{Ind}(v_\a^{(r)},V_\a;p)$ of $v_\a^{(r)}$ on $\sigma\cap V_\a$, if $\dim V_\a \ge r$.
\enddefinition

We now give the proposed proof of the Proportionality Theorem in [BS] in the case of
frames.

\proclaim{Theorem 2.3} Let $V_\a\subset V$ be a stratum and $v_\a^{(r)}$ an $r$-field on 
$\sigma\cap V_\a$ with an isolated 
singularity at the barycenter $p$ of $\sigma$, where $\sigma$ is a cell of dimension 
$2(m-r+1)$. Let $v^{(r)}$ denote an $r$-field on $\sigma$ obtained from $v_\a^{(r)}$ by
radial extension. Then we have
$$
\op{Eu}(v^{(r)},V;p)=\op{Eu}(V,p)\cdot \op{Ind}_{\op{Sch}}(v^{(r)},V;p).
$$
\endproclaim

\demo{Proof}
If $r=1$, this is Theorem 1.3.  If $r>1$, we reduce the problem to the case $r=1$ 
in the following way. 

First, we may write 
$v^{(r)}$ as $(v^{(r-1)},v_r)$, where the $(r-1)$-field $v^{(r-1)}$
is non-singular on $\sigma$. Let $E$ denote the  trivial subbundle of
$TM|_\sigma$ of rank $r-1$ spanned by $v^{(r-1)}$ (over the complex numbers)
and $Q$ the orthogonal 
complement of $E$ in $TM|_\sigma$ for some metric;
$$
TM|_\sigma=E\oplus Q.
$$
Accordingly, we have a decomposition on $\nu^{-1}\sigma$:
$$
\nu^*TM|_\sigma=\nu^*E\oplus \nu^*Q.
$$ 

Since the $r$-field  $(v^{(r-1)},v_r)$ is stratified and non-singular on
$\partial\sigma\cap V$, it lifts to an $r$-frame
$\widetilde v^{(r)} = (\widetilde v^{(r-1)}, \widetilde v_r)$ of the Nash bundle 
$\widetilde T$ over $\nu^{-1}(\partial\sigma\cap V)$. Moreover, since $v^{(r-1)}$
is non-singular on $\sigma\cap V$, $\nu^* E$ (restricted to $\nu^{-1}(\sigma\cap V)$)
is a subbundle of $\widetilde T|_{\nu^{-1}(\sigma\cap V)}$ and we have  a decomposition: 
$$
\widetilde T|_{\nu^{-1}(\sigma\cap V)} = 
\nu^* E \oplus \widetilde P,\tag2.4
$$ 
where $\widetilde P$ is a subbundle of $\nu^*Q|_{\nu^{-1}(\sigma\cap V)}$.
We may think of $\widetilde v_r$ as a section of $\widetilde P$ which is non-vanishing
on $\nu^{-1}(\partial\sigma\cap V)$. If we denote by $o(\widetilde v_r,\widetilde P)$
the class in  
$H^{2(n-r+1)}(\nu^{-1}(\sigma\cap V),\nu^{-1}(\partial\sigma\cap V))$
of the obstruction cocycle to extending this to a non-vanishing section over 
$\nu^{-1}(\sigma\cap V)$, we have
$$
\op{Eu}(v^{(r)},V;p)=o(\widetilde v_r,\widetilde P)
[\nu^{-1}(\sigma\cap V),\nu^{-1}(\partial\sigma\cap V)].\tag2.5
$$

Since the vector field $v_r$ is defined only on $\sigma$, we complement it
by a radial vector field on the complementary space, in order to apply Theorem 1.3.
Thus let $\Bbb D$ denote a small
closed disc of complex dimension $r-1$ with center $p$ and transverse to $\sigma$ and set 
$\Bbb B=\Bbb D\times\sigma$.
The bundle $E$  extends to a neighborhood 
of $p$ in $M$ and we get, denoting also by $E$ and $Q$ the extensions of $E$ and $Q$, 
a decomposition
$$
TM|_{\Bbb B}=E\oplus Q.
$$

The bundles $E$ and $Q$ can be interpreted 
as $\pi_1^*T\Bbb D$ and $\pi_2^*T\sigma$, respectively, where
$\pi_1:\Bbb B\to\Bbb D$ and $\pi_2:\Bbb B\to \sigma$ denote the projections.
Since $E$ is in $SV$, i.e., the vectors in $E$ are stratified, the pull-back 
$\nu^*E$ (restricted to $\nu^{-1}(\Bbb B\cap V))$ is a
subbundle of $\widetilde T|_{\nu^{-1}(\Bbb B\cap V)}$ and we have a decomposition
$$
\widetilde T|_{\nu^{-1}(\Bbb B\cap V)} = 
\nu^* E \oplus \widetilde P,
$$ 
where $\widetilde P$ is a subbundle of $\nu^*Q|_{\nu^{-1}(\Bbb B\cap V)}$, 
extending $\widetilde P$ in (2.4).

Now we may think of $v_r$ as a section of $T\sigma$.
Let $v_{\Bbb D}$ denote a radial vector field on $\Bbb D$ at $p$.
Then the direct sum $v_0 = (\pi_1^*v_{\Bbb D}, \pi_2^*v_r)$ is a stratified  
vector field at $p$. Let 
$\widetilde v_0 = (\widetilde {\pi_1^*v_{\Bbb D}}, \widetilde {\pi_2^*v_r})$
be the lifting of these vector fields to sections of $\widetilde T$ over
$\nu^{-1}(\Bbb S \cap V)$, where $\widetilde {\pi_1^*v_{\Bbb D}}=\nu^*\pi_1^*v_{\Bbb D}$
is a section of $\nu^* E$ and $\widetilde {\pi_2^*v_r}$ is a section of $\widetilde P$. 
If we denote by $o(\widetilde v_0,\widetilde T)$ the class of  obstruction
cocycle to extending $\widetilde v_0$ to a non-vanishing section over 
$\nu^{-1}(\Bbb B\cap V)$, by definition we have
$$
\op{Eu}(v_0,V;p) = o(\widetilde v_0,\widetilde T)
[\nu^{-1}(\Bbb B\cap V),\nu^{-1}(\Bbb S\cap V)].
$$

We now show the identity
$$
\op{Eu}(v^{(r)},V;p)=\op{Eu}(v_0,V;p).\tag2.6
$$ 

Denoting by $ o(\nu^*\pi_1^*v_{\Bbb D},\nu^*E)$ the class of  obstruction
cocycle to extending $\nu^*\pi_1^*v_{\Bbb D}$, a section of $\nu^*E$ non-vanishing on
$\nu^{-1}((\partial\Bbb D\times\sigma)\cap V)=\nu^{-1}(\partial\Bbb D\times (\sigma\cap V))$, 
to a non-vanishing section on $\nu^{-1}(\Bbb B\cap V)$ and by
$o(\widetilde {\pi_2^*v_r},\widetilde P)$ the class of  obstruction
cocycle to extending $\widetilde {\pi_2^*v_r}$, a section of $\widetilde P$ non-vanishing on
$\nu^{-1}((\Bbb D\times\partial\sigma)\cap V)=\nu^{-1}(\Bbb D\times (\partial\sigma\cap V))$, 
to a non-vanishing section on $\nu^{-1}(\Bbb B\cap V)$, we have
$$
 o(\widetilde v_0,\widetilde T)= o(\nu^*\pi_1^*v_{\Bbb D},\nu^*E)\smile  
o(\widetilde {\pi_2^*v_r},\widetilde P),
$$
where $\smile$ denotes the cup product.
We have $o(\nu^*\pi_1^*v_{\Bbb D},\nu^*E)=\nu^*\pi_1^*o(v_{\Bbb D},T\Bbb D)$.
Since $\Bbb B\cap V=\Bbb D\times(\sigma\cap V)$ and $o(v_{\Bbb D},T\Bbb D)$ is a
generator of $H^{2r-2}(\Bbb D,\partial\Bbb D)$, we get
$$
o(\nu^*\pi_1^*v_{\Bbb D},\nu^*E)\frown [\nu^{-1}(\Bbb B\cap V),\nu^{-1}(\Bbb S\cap V)]=
[\nu^{-1}(\sigma\cap V),\nu^{-1}(\partial\sigma\cap V)],
$$
where $\frown$ denotes the cap product.
Since the restriction of $o(\widetilde {\pi_2^*v_r},\widetilde P)$ to $\nu^{-1}\sigma$
is equal to $o(\widetilde {v_r},\widetilde P)$, comparing with (2.5), we obtain (2.6).

By Theorem 1.3, we have $\op{Eu}(v_0,V;p)=\op{Eu}(V,p)\cdot \op{Ind}_{\op{Sch}}(v_0,V;p)$.
Also from definition, we have $\op{Ind}_{\op{Sch}}(v_0,V;p)=\op{Ind}_{\op{Sch}}(v_r,V;p)
=\op{Ind}_{\op{Sch}}(v^{(r)},V;p)$
and the theorem.
\qed
\enddemo

\enddocument
\end